\documentclass{article}

\usepackage{proof}
\usepackage{latexsym}

\newtheorem{theorem}{Theorem}[section]
\newtheorem{definition}{Definition}[section]
\newtheorem{proposition}{Proposition}[section]
\newtheorem{lemma}{Lemma}[section]

\newtheorem{claim}{Claim}[section]

\title{Intuitionistic fixed point theories over Heyting arithmetic
\thanks{Dedicated to the occasion of Grisha Mints' 70th birthday}
}
\author{
Toshiyasu Arai
\\
Toshiyasu ARAI
\\
Graduate School of Science
\\
Chiba University
\\
1-33, Yayoi-cho, Inage-ku,
Chiba, 263-8522, JAPAN
}
\date{}

\begin{document}
\maketitle
\begin{abstract}
In this paper we show that an intuitionistic theory for fixed points is conservative over the Heyting arithmetic with respect to a certain class of formulas.
This extends partly the result of mine.
The proof is inspired by the quick cut-elimination due to G. Mints
\end{abstract}

\section{Introduction}

Fixed points occur frequently in mathematical reasonings.
Let us consider in this paper the fixed point predicate $I^{\Phi}(x)$ for positive formula $\Phi(X,x)$:
\begin{equation}\label{eq:fix}
(FP)^{\Phi} \; \; \forall x[I^{\Phi}(x)  \leftrightarrow  \Phi(I^{\Phi},x)]
\end{equation}

Over the classical logic, the existence of fixed points strengthens theories.

In \cite{Motohashi} we have shown that a first-order logic calculus with the axioms (\ref{eq:fix})
has non-elementary speed-ups over the classical first-order predicate logic.

As an extension of the first-order arithmetic PA,  the theory $\widehat{ID}$ for fixed points is stronger than PA.
In $\widehat{ID}$ one can readily define the truth definition of arithmetic formulas.
Moreover $\widehat{ID}$ proves the transfinite induction up to each ordinal less than $\varphi\varepsilon_{0}0$ for arithmetic formulas, \cite{Feferman} and \cite{attic}.

However {\it intuitionistic\/} theories for fixed points may be proof-theoretically equivalent to 
the intuitionistic arithmetic HA.

W. Buchholz\cite{Buchholz} showed that an intuitionistic fixed point theory $\widehat{ID}^{i}({\cal M})$ is
conservative over the Heyting arithmetic HA with respect to almost negative formulas(, in which
$\lor$ does not occur and $\exists$ occurs in front of atomic formulas only).
The theory $\widehat{ID}^{i}({\cal M})$ has the axioms (\ref{eq:fix}) $(FP)^{\Phi}$ for fixed points
for {\it monotone formula\/} $\Phi(X,x)$, which is generated from arithmetic atomic formulas and $X(t)$ by means of (first order) monotonic connectives $\lor,\land,\exists,\forall$.
Namely $\to$ nor $\lnot$ does occur in monotone formula.
The proof is based on a recursive realizability interpretation.

After seeing the result of Buchholz, we\cite{attic} showed that
an intuitionistic fixed point (second order) theory is conservative over HA for any arithmetic formulas.
In the theory the operator $\Phi$ for fixed points is generated from  $X(t)$ and 
any second order formulas by means of first order monotonic connectives and second order 
existential quantifiers $\exists f(\in\omega\to\omega)$.
The proof in \cite{attic} is to interpret the fixed points by $\Sigma^{1}_{1}$-formulas as in \cite{Feferman}.
In interpreting the fixed points by $\Sigma^{1}_{1}$-formulas, we need an axiom of choice 
$\mbox{AC}_{01}$.
Therefore the proof does not work for strictly positive operators, e.g.,
$\Phi(X,x): \Leftrightarrow \lnot\exists y\forall z A(x) \to X(x)$
since 
$\lnot\exists y\forall z A(x) \to \exists f R(f,x)\leftrightarrow \exists f[\lnot\exists y\forall z A(x)\to R(x)]$
is nothing but the independence of premiss, IP, which is not valid intuitionistically.

The crux is the fact that the axiom of choice adds nothing to HA, i.e., N. Goodman's theorem\cite{Goodman},
while the theorem is proved by a combination of a realizability interpretation and a forcing.
Also cf. \cite{Mintsfinite} for a proof-theoretic proof of the Goodman's theorem.

I met Grisha first time in Hiroshima, Japan, September, 1995.
I explained to him the result in \cite{attic}.
He soon realized that it be followed from the Goodman's theorem
before I told the proof.
Then he asked me "Can you prove it by means of proof-transformations, e.g., cut-elimination?".
This paper is a partial answer to his question.

Now let $\widehat{ID}^{i}({\cal HM})$ denote an intuitionistic fixed point theory in which the operator 
$\Phi(X,x)$ is in a class ${\cal HM}$ of formulas, cf. Definition \ref{df:leftright}  below.
The class ${\cal HM}$ contains properly the monotone formulas
and typically is of the form $H(x) \to M(X,x)$ for a (Rasiowa-)Harrop formula $H$(, in which there
is no strictly positive occurrence of disjunction nor existential subformulas) and a monotonic formula $M$.

We show that the theory $\widehat{ID}^{i}({\cal HM})$ is conservative over HA with respect to
the class ${\cal HM}$.
Thus the result of the paper extends partly one in \cite{attic}.

On the other side, C. R\"uede and T. Strahm\cite{Strahm} extends significantly 
the results in \cite{Buchholz} and \cite{attic}.
They showed that the intuitionistic fixed point theory for {\it strictly positive\/} operators
is conservative over HA with respect to negative and $\Pi^{0}_{2}$-formulas.
Moreover they determined the proof-theoretic strengths of intuitionisitic theories for
transfinitely iterations of fixed points by strictly positive operators.
The class of strictly positive formulas is wider than our class ${\cal HM}$.
In this respect the result in \cite{Strahm} supersedes ours.
A merit here is that the class ${\cal HM}$ is wider than the class concerned in \cite{Strahm}.
For example any formula in prenex normal form is (equivalent to a formula) in ${\cal HM}$, but
a $\Pi^{0}_{3}$-formula is neither negative nor $\Pi^{0}_{2}$.

Rather, I think that the novelty lies in our {\it proof technique\/},
which shows that, cf. Theorem \ref{lem:quickelim},
eliminating cut inferences with ${\cal HM}$-cut formulas from derivations of ${\cal HM}$-end formulas
blows up depths of derivations {\it only by one exponential\/}, e.g., towers of exponentials are dispensable.
This is seen from the fact that there exists an embedding from the resulting tree of cut-free derivation
 to such a derivation with cut inferences such that the embedding maps the deeper nodes in the tree ordering
 to larger nodes with respect to {\it Kleene-Brouwer ordering\/}.
 In other words eliminating monotone cut formulas is essentially to linearize the well founded tree
as in Kleene-Brouwer ordering.
 This is an essence of quick cut-elimination in \cite{Mintsmono}.
 
Let us explain an idea of our proof more closely.
First the finitary derivations in $\widehat{ID}^{i}({\cal HM})$ is embedded to infinitary derivations,
and eliminate cuts partially.
This results in an infinitary derivation of depth less than $\varepsilon_{0}$, and in which there
occurs cut inferences with cut formulas $I^{\Phi}(t)$ for fixed points only.
Now the constrains on operator $\Phi$ and the end formula admits us to
invert cut-free derivations of sequents with a Harrop antecedent and a monotonic succedent formula.
Therefore the quick cut-elimination (and pruning) technique in Grisha's\cite{Mintsmono}
could work to eliminate cut inferences with cut formulas $I^{\Phi}(t)$.
In this way we will get an infinitary derivation of depth less than $\varepsilon_{0}$, and in which there
occurs no fixed point formulas.

By formalizing the arguments we see that the end formula is true in HA.

\section{An intuitionistic theory $\widehat{ID}^{i}({\cal HM})$}

$L_{HA}$ denote the language of the Heyting arithmetic.
$L_{HA}$ consists of the equality sign $=$, individual constants $0,1$ for zero and one, 
function symbols $+,\cdot$ for addition and multiplication, and logical connectives $\lor,\land,\to, \exists,\forall$.

Let $X$ be a fresh predicate symbol, which is assumed to be unary for simplicity.
$L_{HA}(X)$ denotes the language $L_{HA}\cup\{X\}$.

\begin{definition}\label{df:leftright}
{\rm Define inductively two classes of formulas} ${\cal H}$ {\rm in} $L_{HA}$,
 {\rm and} ${\cal HM}$ {\rm in} $L_{HA}(X)$
{\rm as follows.}
\begin{enumerate}
\item
{\rm Any atomic formula} $s=t$ {\rm belongs to both of} ${\cal H}$ {\rm and} ${\cal HM}$.
\item
{\rm Any atomic formula} $X(t)$ {\rm belongs to the class} ${\cal HM}$.
\item
{\rm If} $H,G\in{\cal H}${\rm , then} $H\land G, \forall x H\in{\cal H}$.
\item
{\rm If} $H\in{\cal H}${\rm , then} $A\to H\in{\cal H}$ {\rm for any formula} $A\in L_{HA}$.
\item
{\rm If} $R,S\in{\cal HM}${\rm , then} $R\lor S, R\land S, \exists x R,\forall x R\in{\cal HM}$.
\item\label{df:leftright9}
{\rm If} $L\in{\cal H}$ {\rm and} $R\in{\cal HM}${\rm , then} $L\to R\in{\cal HM}$.
\end{enumerate}
\end{definition}

${\cal H}$ denotes the class of (Rasiowa-)Harrop formulas, in which there occurs no strictly positive existential nor disjunctive subformula.

${\cal HM}$ contains properly the monotone formulas, i.e., the class POS in \cite{Buchholz},
e.g., $\lnot A\to X(a)\in{\cal HM}$ is not intuitionistically equivalent to any monotone formula,
but there exists a strongly positive formula with respect to $X$ not in ${\cal HM}$, e.g.,
$(\forall x\exists y A\to\exists z B)\land X(a)$.
Any formula in ${\cal HM}$ is strictly positive with respect to $X$.

Let $\widehat{ID}^{i}({\cal HM})$ denote the following extension of HA.
Its language is obtained from $L_{HA}$ by adding a unary set constant $I^{\Phi}$ for each 
$\Phi\equiv\Phi(X,x)\in {\cal HM}$, in which only a fixed variable $x$ occurs freely.
Its axioms are those of HA in the expanded language(, i.e., the induction axioms are available for
any formulas in the expanded language)  plus the axiom $(FP)^{\Phi}$, (\ref{eq:fix}) for fixed points.

Now our theorem runs as follows.

\begin{theorem}\label{th:main}
$\widehat{ID}^{i}({\cal HM})$ is conservative over {\rm HA} with respect to formulas in ${\cal HM}$
(, in which the extra predicate constant $X$ does not occur).
\end{theorem}

\section{Infinitary derivations}

Given an $\widehat{ID}^{i}({\cal HM})$-derivation $D_{0}$ of an ${\cal HM}$-sentence $R_{0}$, let us first embed it to
an infinitary derivation in an infinitary calculus $\widehat{ID}^{i\infty}({\cal HM})$.

\[\lnot A:\Leftrightarrow A\to \bot.\]

Let $N$ denote a number which is big enough so that any formula occurring in $D_{0}$ has
logical complexity(, which is defined by the number of occurrences of logical connectives) smaller than $N$.
In what follows any formula occurring in infinitary derivations which we are concerned, has
logical complexity less than $N$.

The derived objects in the calculus $\widehat{ID}^{i\infty}({\cal HM})$ are {\it sequents\/}
$\Gamma\Rightarrow A$, where $A$ is a {\it sentence\/} (in the language of $\widehat{ID}^{i}({\cal HM})$)
and $\Gamma$ denotes a finite set of {\it sentences\/},
where each closed term $t$ is identified with its value $\bar{n}$, the $n$th numeral.

$\bot$ stands ambiguously for false equations $t=s$ with  closed terms $t,s$ having different values.
$\top$ stands ambiguously for true equations $t=s$ with closed terms $t,s$ having same values.

The {\it initial sequents\/} are
\[
\Gamma,I(t)\Rightarrow I(t)\, ; \mbox{\hspace{5mm}} \Gamma,\bot\Rightarrow A\, ; \Gamma\Rightarrow \top
\]
These are regarded as inference rules with empty premiss(upper sequent).

The {\it inference rules\/} are $(L\lor)$, $(R\lor)$, $(L\land)$, $(R\land)$, $(L\to)$, $(R\to)$, 
$(L\exists)$, $(R\exists)$, $(L\forall)$, $(R\forall)$, $(LI)$, $(RI)$, $(cut)$, and the repetition rule $(Rep)$.
These are standard ones.

\begin{enumerate}

\item
\[
\infer[(LI)]
{\Gamma,I(t) \Rightarrow C}
{\Gamma,\Phi(I,t) \Rightarrow C}
\: ;\:
\infer[(RI)]
{\Gamma\Rightarrow I(t)}
{\Gamma\Rightarrow \Phi(I,t)}
\]
\item
\[
\infer[(L\lor)]
{\Gamma,A_{0}\lor A_{1}\Rightarrow C}
{
\Gamma,A_{0}\Rightarrow C
&
\Gamma,A_{1}\Rightarrow C
}
\: ;\:
\infer[(R\lor)]
{\Gamma\Rightarrow A_{0}\lor A_{1}}
{\Gamma\Rightarrow A_{i}}
\,(i=0,1)
\]
\item
\[
\infer[(L\land)]
{\Gamma,A_{0}\land A_{1}\Rightarrow C}
{
\Gamma,A_{0}\land A_{1},A_{i}\Rightarrow C
}
\, (i=0,1)
\: ;\:
\infer[(R\land)]
{\Gamma\Rightarrow A_{0}\land A_{1}}
{
\Gamma\Rightarrow A_{0}
&
\Gamma\Rightarrow A_{1}
}
\]
\item
\[
\infer[(L\to)]
{\Gamma,A\to B\Rightarrow C}
{
\Gamma,A\to B\Rightarrow A
&
\Gamma,B\Rightarrow C
}
\: ;\:
\infer[(R\to)]
{\Gamma\Rightarrow A\to B}
{\Gamma,A\Rightarrow B}
\]
 \item
\[
\infer[(L\exists)]
{\Gamma,\exists x B(x)\Rightarrow C}
{
\cdots
&
\Gamma,B(\bar{n})\Rightarrow C
&
\cdots (n\in\omega)
}
\: ;\:
\infer[(R\exists)]
{\Gamma\Rightarrow \exists x B(x)}
{\Gamma\Rightarrow B(\bar{n})}
\]
 \item
 \[
 \infer[(L\forall)]
 {\Gamma,\forall x B(x)\Rightarrow C}
 {\Gamma,\forall x B(x),B(\bar{n})\Rightarrow C}
 \: ;\:
 \infer[(R\forall)]
 {\Gamma\Rightarrow \forall x B(x)}
 {
 \cdots
 &
 \Gamma\Rightarrow B(\bar{n})
 &
 \cdots (n\in\omega)
 }
 \]
  \item
  \[
  \infer[(cut)]
  {\Gamma,\Delta\Rightarrow C}
  {
  \Gamma\Rightarrow A
  &
  \Delta,A\Rightarrow C
  }
  \]
  \item
  \[
  \infer[(Rep)]
  {\Gamma\Rightarrow C}{\Gamma\Rightarrow C}
  \]
\end{enumerate}

The {\it depth\/} of an infinitary derivation is defined to be the depth of the well founded tree.

As usual we see the following proposition.
Recall that $N$ is an upper bound of logical complexities of formulas occurring in the given finite derivation $D_{0}$ of ${\cal HM}$-sentence $R_{0}$.

\begin{proposition}\label{prp:embed}
\begin{enumerate}
\item\label{prp:embed1}
There exists an infinitary derivation $D_{1}$ of $R_{0}$ such that its depth is less than $\omega^{2}$ and
the logical complexity of any sentence, in particular cut formulas occurring in $D_{1}$ is less than $N$.
\item\label{prp:embed2}
By a partial cut-elimination, there exist an infinitary derivation $D_{2}$ of $R_{0}$ and
an ordinal $\alpha_{0}<\varepsilon_{0}$
such that
the depth of the derivation $D_{2}$ is less than $\alpha_{0}$ and any cut formula occurring in $D_{2}$ is an atomic formula
$I(t)$(, and the logical complexity of any formula occurring in it is less than $N$).
\end{enumerate}
\end{proposition}

Let ${\cal HM}(I)$ denote the class of sentences obtained from sentences in ${\cal HM}$ by substituting
the predicate $I^{\Phi}$ for the predicate $X$.
\begin{definition}\label{df:rank}
{\rm The} rank $rk(A)$ {\rm of a sentence} $A$ {\rm is defined by}
\[
rk(A):= \left\{
\begin{array}{ll}
0 & \mbox{{\rm if} } A\in{\cal H} \\
1 & \mbox{{\rm if }} A\in {\cal HM}(I)\setminus{\cal H} \\
2 & \mbox{{\rm otherwise}}
\end{array}
\right.
\]

{\rm For inference rules} $J${\rm , the} rank $rk(J)$ {\rm of} $J$ 
{\rm is defined to be the rank of the cut formula if} $J$ {\rm is a cut inference.
Otherwise} $rk(J):=0$.

{\rm For derivations} $D${\rm , the} rank $rk(D)$ {\rm of} $D$ 
{\rm is defined to be the maximum rank of the cut inferences in it.}
\end{definition}

Let $\vdash^{\alpha}_{r}\Gamma\Rightarrow C$ mean that there exists an infinitary derivation of $\Gamma\Rightarrow C$
such that its depth is at most $\alpha$, and its rank is less than $r$(, and
and the logical complexity of any formula occurring in it is less than $N$).

\begin{theorem}\label{lem:quickelim}
Let $C_{0}$ denote an ${\cal HM}$-sentence, and $\Gamma_{0}$ a finite set of ${\cal H}$-sentences.
Suppose that $\vdash^{\alpha}_{2}\Gamma_{0}\Rightarrow C_{0}$.
Then $\vdash^{\omega^{\alpha}+1}_{1}\Gamma_{0}\Rightarrow C_{0}$.
\end{theorem}

Assuming the Theorem \ref{lem:quickelim}, we can show the Theorem \ref{th:main} as follows.
Suppose an ${\cal HM}$-sentence $C_{0}$ is provable in $\widehat{ID}^{i}({\cal HM})$.
By Proposition \ref{prp:embed} we have $\vdash^{\alpha_{0}}_{2}\Rightarrow C_{0}$ for a big enough 
number $N$ and an $\alpha_{0}<\varepsilon_{0}$.
Then Theorem \ref{lem:quickelim} yields $\vdash^{\beta_{0}}_{1}\Rightarrow C_{0}$
for $\beta_{0}=\omega^{\alpha_{0}}+1<\varepsilon_{0}$.

Let ${\rm Tr}_{N}(x)$ denote a partial truth definition for formulas of logical complexity less than $N$,
cf. \cite{Troelstra}, 1.5.4.
By transfinite induction up to $\beta_{0}$, cf. Lemma \ref{lem:KB}, we see ${\rm Tr}_{N}(C_{0})$.
Note that any sentence occurring  in the witnessed derivation for $\vdash^{\beta_{0}}_{1}\Rightarrow C_{0}$
has logical complexity less than $N$, and it is either an ${\cal H}$-sentence or an ${\cal HM}$-sentence.
Specifically there occurs no fixed point formula $I(t)$ in it.
Now since everything up to this point is formalizable in HA, we have ${\rm Tr}_{N}(C_{0})$, and hence $C_{0}$ in HA.
This shows the Theorem \ref{th:main}.

A proof of Theorem \ref{lem:quickelim} is given in the next section.

\section{Quick cut-elimination for monotone cuts with Harrop side formulas}

Our plan of the proof of Theorem \ref{lem:quickelim} is as follows.
Pick the leftmost cut $J$ of rank 1:
\[
\infer[(cut)J]
{\Gamma,\Delta\Rightarrow C}
{
 \infer*[D_{\ell}]{\Gamma\Rightarrow A}{}
& 
 \infer*[D_{r}]{\Delta,A\Rightarrow C}{}
}
\]
with $rk(A)=1$.

Then $\vdash^{\alpha}_{1}\Gamma\Rightarrow A$ and $\vdash^{\beta}_{2}\Delta,A\Rightarrow C$ for some $\alpha$ and $\beta$.
Moreover since for the end sequent $\Gamma_{0}\Rightarrow C_{0}$, $\Gamma_{0}\subseteq{\cal H}$ and $C_{0}\in{\cal HM}$,
we see that any sentence in $\Gamma\cup\Delta$ is in ${\cal H}$, and $C\in{\cal HM}$ by 
$\vdash^{\alpha_{0}}_{2}\Gamma_{0}\Rightarrow C_{0}$.

Since $\Gamma$ consists solely in Harrop formulas,
we can invert the derivation $D_{\ell}$, and climb up the derivation $D_{r}$ with inverted $D_{\ell}$.
This results in a derivation of $\Gamma,\Delta\Rightarrow C$ of depth $dp(D_{\ell})+dp(D_{r})$.
Iterating this eliminations, we could get a derivation of rank 0, and of depth at most exponential of the depth of the given derivation.

Though intuitively this would suffice to believe in Theorem \ref{lem:quickelim},
we have to prove two facts:
first why the iteration eventually terminates?
second give a succinct argument to the estimated increase of depth.
These are not entirely trivial tasks.
It turns out that we need to proceed along the Kleene-Brouwer ordering on well founded trees
 instead of the depths.
Let us explore this.

Let us fix a witnessed derivation $D_{2}$ of $\vdash^{\alpha_{0}}_{2}\Gamma_{0}\Rightarrow C_{0}$.
Let $(T_{2},<_{T_{2}})$ denote the wellordering, where $T_{2}\subseteq{}^{<\omega}\omega$ is
the naked tree of $D_{2}$ and $<_{T_{2}}$ the Kleene-Brouwer ordering on $T_{2}$.

Let us consider infinitary derivations equipped with additional informations as in \cite{Mintsfinite}.
\begin{definition}\label{df:derivation}
{\rm An} infinitary derivation {\rm is a sextuple} $D=(T,Seq,Rule, rk,ord,kb)$ {\rm which enjoys the following conditions.
The naked tree of} $D$ {\rm is denoted} $T=T(D)$.
\begin{enumerate}
\item
$T\subseteq{}^{<\omega}\omega$ {\rm is a tree in the sense that
there exists a root} $r\in T$ {\rm with} $\forall a\in T(r\subseteq a)$ {\rm and}
$\forall a, b(r\subseteq a\subseteq b\in T \Rightarrow a\in T)$.

{\rm It is} not assumed {\rm that the empty node} $\emptyset$ {\rm is to be the root nor}
$a*\langle n\rangle\in T \,\&\, m<n \Rightarrow a*\langle m\rangle\in T$.

\item
$Seq(a)$ {\rm for} $a\in T$ {\rm denotes the sequent situated at the node} $a$.

{\rm If} $Seq(a)$ {\rm is a sequent} $\Gamma\Rightarrow C${\rm , then it is denoted}
\[
a:\Gamma\Rightarrow C
.\]

\item
$Rule(a)$ {\rm for} $a\in T$ {\rm denotes the name of the inference rule with its lower sequent} $Seq(a)$.
\item
$rk(a)$ {\rm for} $a\in T$ {\rm denotes the rank of the inference rule} $rk(a)$.
\item
$ord(a)$ {\rm for} $a\in T$ {\rm denotes the ordinal}$<\varepsilon_{0}$ {\rm attached to} $a$.
\item
{\rm The quintuple} $(T,Seq,Rule, rk,ord)$ {\rm has to be locally correct with respect to}
 $\widehat{ID}^{i\infty}({\cal HM})$ {\rm and
for being well founded tree} $T$.
\end{enumerate}
{\rm Besides these conditions an extra information is provided by a} labeling function $kb$.

$kb:T\to T_{2}$ {\rm is a function such that for} $a,b\in T$
\begin{enumerate}
\item
\begin{equation}\label{eq:kb1}
kb(a)\neq kb(b) \Rightarrow [a<_{T}b \Leftrightarrow kb(a)<_{T_{2}}kb(b)]
\end{equation}
{\rm for the} Kleene-Brouwer ordering $<_{T}$ {\rm on} $T$.
\item
\begin{equation}\label{eq:kb2}
a\subset b \Rightarrow kb(b)<_{T_{2}}kb(a)
\end{equation}
{\rm where} $a\subset b$ {\rm means that} $a$ {\rm is a proper initial segment of} $b$ 
{\rm for} $a,b\in{}^{<\omega}\omega$.
\item
{\rm Let} $c\in T$ {\rm be a node with} $rk(c)=1${\rm , and} $a,b\in T$ {\rm nodes such that}
$c*\langle \ell\rangle\subseteq a$ {\rm and} $c*\langle r\rangle\subseteq b$ {\rm for} $\ell<r$.
{\rm (This means that} $Seq(a)$ {\rm [}$Seq(b)${\rm ] is in the left [right] upper part of the cut inference}
$Rule(c)${\rm .)
Suppose that the} right cut formula $A$ {\rm in the antecedent of} $Seq(c*\langle r\rangle)$ {\rm has an} ancestor {\rm in} 
$Seq(b)$.

{\rm Then}
\begin{equation}\label{eq:kb3}
 kb(a)\neq kb(b)
\end{equation}
\end{enumerate}
\end{definition}
The condition (\ref{eq:kb2}) to $kb$ ensures us that the depth of $T$ is at most the order type of $<_{T_{2}}$.

It is easy to see that Kleene-Brouwer ordering $<_{T_{2}}$ is a well ordering, and
its order type is bounded by $\omega^{\alpha}+1$ for the depth $\alpha$ of the primitive recursive and
wellfounded tree $T_{2}$.

\begin{lemma}\label{lem:KB}
 The transfinite induction schema (for arithmetical formulas) along the Kleene-Brouwer ordering
 $<_{T_{2}}$
 is provable in HA.
\end{lemma}
 {\bf Proof}.\hspace{2mm}
 Since the transfinite induction schema along a standard $\varepsilon_{0}$-ordering is provable in HA up to
 each $\alpha<\varepsilon_{0}$, the same holds for the tree ordering $\{(b,a): a\subset b, a,b\in T_{2}\}$ (bar induction).
 
 Now let $X$ be a formula, and assume that $X$ is progressive with respect to $<_{T_{2}}$:
 \[
 \forall a\in T_{2}[\forall b<_{T_{2}}a\, X(b) \to X(a)]
 .\]
 Let
 \[
 {\sf j}[X](a):\Leftrightarrow \forall y\in T_{2}[y\supseteq a \to \forall x<_{T_{2}}y\, X(x) \to \forall x<_{T_{2}}a\, X(x)]
 .\]
 Then we see that ${\sf j}[X]$ is progressive with respect to the tree ordering.
 Therefore ${\sf j}[X](r)$ for the root $r\in T_{2}$, and by letting $y$ to be the leftmost leaf in $T_{2}$
 we have $\forall x<_{T_{2}}r\, X(x)$. The progressiveness of $X$ with respect to $<_{T_{2}}$
 yields $X(r)$.
 \hspace*{\fill} $\Box$

The following Lemmas are seen as usual.

\begin{lemma}\label{lem:weakfalse}
Let $D=(T,Seq,Rule, rk,ord,kb)$ be a derivation of rank 1, and of a sequent $\Gamma\Rightarrow A$
such that $\Gamma\subseteq{\cal H}$ and $A\in{\cal HM}\,\&\, rk(A)=1$.
For any $\Delta$, there exists a derivation $D*\Delta=(T,Seq*\Delta,Rule, rk,ord,kb)$ of the sequent
$\Gamma,\Delta\Rightarrow A$.
\end{lemma}

\begin{lemma}\label{lem:inversion}(Inversion Lemma)\\
Let $D=(T,Seq,Rule, rk,ord,kb)$ be a derivation of rank 1, and of a sequent $\Gamma\Rightarrow A$
such that $\Gamma\subseteq{\cal H}$ and $A\in{\cal HM}\,\&\, rk(A)=1$.
\begin{enumerate}
\item\label{lem:inversion1}
If $A\equiv B_{0}\lor B_{1}$, then there exists a derivation $D_{i}=(T,Seq_{i},Rule_{i}, rk,ord,kb)$
of rank 1 and of a sequent $\Gamma\Rightarrow B_{i}$ for an $i=0,1$.
\item\label{lem:inversion2}
If $A\equiv B_{0}\land B_{1}$, then there exist derivations $D_{i}=(T_{i},Seq_{i},Rule_{i}, rk,ord,kb)$
of rank 1 and of sequents $\Gamma\Rightarrow B_{i}$ for any $i=0,1$, where
$T_{i}\subseteq T$ by pruning.
\item\label{lem:inversion3}
If $A\equiv \exists x B(x)$, then there exists a derivation $D_{n}=(T,Seq_{n},Rule_{n}, rk,ord,kb)$
of rank 1 and of a sequent $\Gamma\Rightarrow B(\bar{n})$ for an $n\in\omega$.
\item\label{lem:inversion4}
If $A\equiv \forall x B(x)$, then there exist derivations $D_{n}=(T_{n},Seq_{n},Rule_{n}, rk,ord,kb)$
of rank 1 and of a sequent $\Gamma\Rightarrow B(\bar{n})$ for any $n\in\omega$, where
$T_{n}\subseteq T$ by pruning.
\item\label{lem:inversion5}
If $A\equiv B_{0}\to B_{1}$, then there exist a derivation $D^{\prime}=(T,Seq^{\prime},Rule^{\prime}, rk,ord,kb)$
of rank 1 and of sequents $\Gamma,B_{0}\Rightarrow B_{1}$.
\item\label{lem:inversion6}
If $A\equiv I(t)$, then there exists a derivation $D^{\prime}=(T,Seq^{\prime},Rule^{\prime}, rk,ord,kb)$
of rank 1 and of sequents $\Gamma\Rightarrow \Phi(I,t)$.
\end{enumerate}
\end{lemma}

\begin{definition}\label{df:KBJ}
{\rm For each cut inference} $J$ {\rm in a derivation} $D$, 
$KB(J)$ {\rm denotes} $kb(a)\in T_{2}$ {\rm with the} left upper node $a$ {\rm of} $J${\rm :}
\[
\infer[J]{\Gamma,\Delta\Rightarrow C}
{
 a :\Gamma\Rightarrow A
 &
 \Delta,A\Rightarrow C
}
\]
\end{definition}

Let us define a cut-eliminating operator $ce_{1}(D)$ for derivations $D=(T,Seq,Rule, rk,ord,kb)$
of rank 1 and of an end sequent $\Gamma_{0}\Rightarrow C_{0}$ with $\Gamma_{0}\subseteq{\cal H}$ and $C_{0}\in {\cal HM}$.

If $D$ is of rank 0, then $ce_{1}(D):=D$.

Assume that $D$ contains a cut inference of rank 1.
Pick the leftmost cut of rank 1:
\[
D=
\infer*{\Gamma_{0}\Rightarrow C_{0}}
{
\infer[(cut)]
{\Gamma,\Delta\Rightarrow C}
{
 \infer*[D_{\ell}]{\Gamma\Rightarrow A}{}
& 
 \infer*[D_{r}]{\Delta,A\Rightarrow C}{}
}
}
\]
The leftmostness means that $KB(J)$ is least in the Kleene-Brouwer ordering $<_{T_{2}}$.

By recursion on the depth
\footnote{As in \cite{Mintsfinite} we see that the operators $ce_{1}, ce_{2}$ are primitive recursive.
We don't need this fact.}
of the derivation $D_{r}$ we define a derivation $ce_{2}(D_{\ell},D_{r})$ of $\Gamma,\Delta\Rightarrow C$.
Then $ce_{1}(D)$ is obtained from $D$ by pruning $D_{\ell}$ and replacing $D_{r}$ by $ce_{2}(D_{\ell},D_{r})$, i.e.,
by grafting $ce_{2}(D_{\ell},D_{r})$ onto the trunk of $D$ up to $\Gamma,\Delta\Rightarrow C$.

As in Lemma 3.2, \cite{Mintsmono} the construction of $ce_{2}(D_{\ell},D_{r})$ is fairly standard,
leaving the resulted cut inferences of rank 0, but has to performed parallely.

Let $\Gamma\cup\Delta\subseteq{\cal H}$ and $\mbox{\boldmath$A$} =A_{1},\ldots,A_{k}$ be a finite sequence of ${\cal HM}$-sentences.
Let $\mbox{\boldmath$D$} _{\ell}=D_{\ell,1},\ldots,D_{\ell,k}$ be rank 0 derivations of
$\Gamma\Rightarrow A_{i}$, and $D_{r}$ a rank 1 derivation of
$\Delta,\mbox{\boldmath$A$} \Rightarrow C$.
We will eliminate the cuts with the cut formulas $A_{i}$ in parallel.
$ce_{2}(\mbox{\boldmath$D$} _{\ell},D_{r})$ is defined from the resulting derivation, denoted $E$ by recursion.

\[
D_{a}=
\infer[(cut)]
{a: \Gamma,\Delta\Rightarrow C}
{
 \infer*[\mbox{\boldmath$D$} _{\ell}]{\mbox{\boldmath$b$} : \Gamma\Rightarrow \mbox{\boldmath$A$} }{}
& 
 \infer*[D_{r}]{c_{1}: \Delta,\mbox{\boldmath$A$} \Rightarrow C}{}
}
\]
denotes the series of cut inferences:
\[
\infer
{a: \Gamma,\Delta\Rightarrow C}
{
 \infer*[D_{\ell,k}]{b_{k}: \Gamma\Rightarrow A_{k}}{}
& 
 \infer*{c_{k}: \Gamma,\Delta,A_{k}\Rightarrow C}
 {
  \infer{c_{2}: \Gamma,\Delta,A_{2},\ldots,A_{k}\Rightarrow C}
  {
   \infer*[D_{\ell,1}]{b_{1}: \Gamma\Rightarrow A_{1}}{}
   &
   \infer*[D_{r}]{c_{1}: \Delta,\mbox{\boldmath$A$} \Rightarrow C}{}
  }
 }
}
\]

\begin{enumerate}
\item\label{1}
If $\Delta,\mbox{\boldmath$A$} \Rightarrow C$ is an initial sequent such that one of the cases $C\equiv\top$, $\bot\in\Delta$ 
or $C\in\Delta$ occurs, then $\Delta\Rightarrow C$, and hence
$\Gamma,\Delta\Rightarrow C$
is still the same kind of initial sequent. For example
\[
\infer[(Rep)]
{a:\Gamma,\Delta\Rightarrow \top}
{
 \infer*{c_{k}: \Gamma,\Delta\Rightarrow \top}
 {
  \infer[(Rep)]{c_{2}: \Gamma,\Delta\Rightarrow \top}
                    {c_{1}: \Gamma,\Delta\Rightarrow \top}
   }
  }
\]
The $T(E)$ is defined by
\[
d\in T(E) \Leftrightarrow d\in T(D) \,\&\, \forall i(b_{i}\not\subseteq d)
.\]

\item\label{2}
If $\Delta,\mbox{\boldmath$A$} \Rightarrow C$ is an initial sequent with the principal formula $\mbox{\boldmath$A$} \ni A_{i}\equiv C\equiv I(t)$,
then $E$ is defined to be
\[
\infer[(Rep)]
{a:\Gamma,\Delta\Rightarrow C}
 {
 \infer*{c_{k}: \Gamma,\Delta\Rightarrow C}
 {
  \infer[(Rep)]{c_{i+1}: \Gamma,\Delta\Rightarrow C}
  {
   \infer*[D_{\ell,i}*\Delta]{b_{i}: \Gamma,\Delta\Rightarrow C}{}
   }
  }
 }
 \]
where $D_{\ell,i}*\Delta$ is obtained from $D_{\ell,i}$ by weakening, cf. Lemma \ref{lem:weakfalse}.

\[
d\in T(E) \Leftrightarrow d\in T(D) \,\&\, \forall j\neq i(b_{j}\not\subseteq d) \,\&\, c_{i}\not\subseteq d
.\]

\item\label{3}
If $A_{i}\in\mbox{\boldmath$A$} $ is of rank 0, then do nothing for the cut inference of $A_{i}$.
 
In each of the above cases $T(E)\subseteq T(D_{a})$.
The labeling function $kb_{E}$ for $E$ is defined to be the restriction of $kb_{D_{a}}$ to $T(E)$.

In what follows assume that $\Delta,\mbox{\boldmath$A$} \Rightarrow C$ is a lower sequent of an inference rule $J$.

\item\label{4}
If the principal formula of $J$ is not in $\mbox{\boldmath$A$} $, then lift up $\mbox{\boldmath$D$} _{\ell}$:
\[
\infer{a:\Gamma,\Delta\Rightarrow C}
{
\mbox{\boldmath$b$} :\Gamma\Rightarrow \mbox{\boldmath$A$} 
&
\infer[(J)]
{c_{1}: \Delta,\mbox{\boldmath$A$} \Rightarrow C}{\cdots & c_{1,i}: \Delta_{i},\mbox{\boldmath$A$} \Rightarrow C_{i} & \cdots}
}
\]
where $\mbox{\boldmath$b$} =b_{k},\ldots,b_{1}$ with $b_{j}=c_{j+1}*\langle \ell_{j}\rangle\, (c_{k+1}:=a)$,
and $c_{j}=c_{j+1}*\langle r_{j}\rangle$ for some $\ell_{j}<r_{j}$, 
and $c_{1,i}=c_{1}*\langle n_{i}\rangle$ for some $n_{i}$ with $i<j\Rightarrow n_{i}<n_{j}$.

$E$ is defined as follows.
\[
\infer[(J)]
{a: \Gamma,\Delta\Rightarrow C}
{
\cdots 
&
 \infer
 {a_{i}:\Gamma,\Delta_{i}\Rightarrow C_{i}}
 {\mbox{\boldmath$b$} _{i}: \Gamma\Rightarrow \mbox{\boldmath$A$}  & c^{\prime}_{1,i}: \Delta_{i},\mbox{\boldmath$A$} \Rightarrow C_{i}}
&
\cdots
}
\]
where $a_{i}=a*\langle n_{i}\rangle$, $\mbox{\boldmath$b$} _{i}=b_{k,i},\ldots,b_{1,i}$ with
$b_{j,i}=c^{\prime}_{j+1,i}*\langle \ell_{j}\rangle\, (c^{\prime}_{k+1,i}:=a_{i})$,
and $c^{\prime}_{j,i}=c^{\prime}_{j+1,i}*\langle r_{j}\rangle$.

The labeling function $kb_{E}$ is defined by
\begin{eqnarray*}
kb_{E}(a*\langle n_{i}\rangle) & = & kb_{D_{a}}(a*\langle r_{k}\rangle)=kb_{D_{a}}(c_{k}),
\\
kb_{E}(a*\langle n_{i}\rangle*\langle r_{k},\ldots, r_{j+1},\ell_{j}\rangle*d)& = & kb_{D_{a}}(a*\langle r_{k},\ldots, r_{j+1},\ell_{j}\rangle*d) \, (1\leq j\leq k),
\\
kb_{E}(a*\langle n_{i}\rangle*\langle r_{k},\ldots, r_{1}\rangle*d)& = & kb_{D_{a}}(a*\langle r_{k},\ldots, r_{1}\rangle*\langle n_{i}\rangle*d)
\end{eqnarray*}

\item\label{5}
Finally suppose that the principal formula of $J$ is a cut formula $A_{i}\in\mbox{\boldmath$A$} $ of $rk(A_{i})=1$.
Use the Inversion Lemma \ref{lem:inversion}.
 \begin{enumerate}
 \item\label{5a}
The case when $A_{i}\equiv \exists x B(x)\in\mbox{\boldmath$A$} $.
For simplicity suppose $i=1$.
\[
\infer[(L\exists)]
{c_{1}: \Delta,\mbox{\boldmath$A$}  \Rightarrow C}
{
 \cdots
 &
 \infer*[D_{r,n}]{c_{1,n}: \Delta, \mbox{\boldmath$A$} _{1}, B(\bar{n})\Rightarrow C}{}
 &
 \cdots
 }
 \]
 where $A_{1}\not\in\mbox{\boldmath$A$} _{1}$.
 
By Inversion Lemma \ref{lem:inversion}.\ref{lem:inversion3} pick an $n$ such that
$\Gamma\Rightarrow B(\bar{n})$ is provable without changing the naked tree.

$E$ is defined as follows.
 \[
 \infer
 {a: \Gamma,\Delta\Rightarrow C}
 {
   \mbox{\boldmath$b$} _{1}: \Gamma\Rightarrow \mbox{\boldmath$A$} _{1}
  &  
  \infer{c_{2}: \Gamma,\Delta,\mbox{\boldmath$A$} _{1}\Rightarrow C}
   {
   b_{1}: \Gamma\Rightarrow B(\bar{n})
   &
   \infer[(Rep)]{c_{1}: \Delta, \mbox{\boldmath$A$} _{1}, B(\bar{n})\Rightarrow C}
   {
   \infer*[D_{r,n}]{c_{1,n}: \Delta, \mbox{\boldmath$A$} _{1}, B(\bar{n})\Rightarrow C}{}
   }
  }
 }
 \]

 \item\label{5b}
 The case when $A_{i}\equiv H\to A_{0}\in\mbox{\boldmath$A$} $ with an $H\in{\cal H}$ and an $A_{0}\in{\cal HM}$.
 For simplicity suppose $i=1$.
\[
\infer[(L\to)]
{c_{1}: \Delta,\mbox{\boldmath$A$}  \Rightarrow C}
{
 c_{1,\ell}: \Delta,\mbox{\boldmath$A$}  \Rightarrow H
&
 c_{1,r}: \Delta,\mbox{\boldmath$A$} _{1},A_{0}\Rightarrow C
}
\]
where $A_{1}\not\in\mbox{\boldmath$A$} _{1}$, and for $m=\ell, r$, 
$c_{1,m}=c_{1}*\langle j_{m}\rangle$ with $j_{\ell}<j_{r}$.

$E$ is defined as follows.
{\footnotesize
\[
\infer
{a: \Gamma,\Delta \Rightarrow C}
 {
  \infer{c_{k,0}: \Gamma,\Delta\Rightarrow H}
  {
  \mbox{\boldmath$b$} _{0}:  \Gamma\Rightarrow \mbox{\boldmath$A$} 
   &
  \infer[(Rep)]{c_{1,0}:  \Delta,\mbox{\boldmath$A$} \Rightarrow H}{c_{1,0,\ell}:  \Delta,\mbox{\boldmath$A$} \Rightarrow H}
   }
  &
  \infer{c_{k,1}: \Gamma,\Delta,H\Rightarrow C}
  {
  \mbox{\boldmath$b$} _{1}: \Gamma\Rightarrow \mbox{\boldmath$A$} _{1}
  &
  \hspace{-10mm}
  \infer{c_{2,1}:  \Gamma,\Delta,H,\mbox{\boldmath$A$} _{1}\Rightarrow C}
   {
   b_{1,1}: \Gamma,H\Rightarrow A_{0}
  &
   \infer[(Rep)]{c_{1,1}:  \Delta,\mbox{\boldmath$A$} _{1}\cup\{A_{0}\}\Rightarrow C}{c_{1,1,r}:  \Delta,\mbox{\boldmath$A$} _{1}\cup\{A_{0}\}\Rightarrow C}
    }
   }
 }
\]
}
where $\Gamma,H\Rightarrow A_{0}$ by inversion.

For $m=0,1$,
$c_{j,m}=c_{j+1,m}*\langle 2r_{j}+m\rangle$ for $1\leq j\leq k$
 with $c_{k+1,m}=a$,
and $\mbox{\boldmath$b$} _{m}=b_{k,m},\ldots,b_{1,m}$ with $b_{j,m}=a*\langle 2r_{k}+m,\ldots, 2r_{j}+m,2\ell_{j}+m\rangle$
and $c_{1,0,\ell}=c_{1,0}*\langle j_{\ell}\rangle$, $c_{1,1,r}=c_{1,1}*\langle j_{r}\rangle$.

The labeling function is defined by
\begin{eqnarray}
kb_{E}(c_{j,m}) & = & kb_{D_{a}}(c_{j})
\label{eq:labelex}
\\
kb_{E}(b_{j,m}*d) & = & kb_{D_{a}}(b_{j}*d), \, (m=0,1)
\nonumber
\end{eqnarray}

\item\label{5c}
 The case when $A_{i}\equiv \forall x B(x) \in\mbox{\boldmath$A$} $.
 For simplicity suppose $i=1$.
 \[
\infer[(L\forall)]
{c_{1}: \Delta,\mbox{\boldmath$A$}  \Rightarrow C}
{
c_{1,n}:  \infer*[D_{r,n}]{\Delta, \mbox{\boldmath$A$} , B(\bar{n})\Rightarrow C}{}
 }
 \]
with $c_{1,n}=c_{1}*\langle j_{n}\rangle$.

$E$ is defined as follows.
 \[
 \infer
 {a: \Gamma,\Delta\Rightarrow C}
 {
   \mbox{\boldmath$b$}  : \Gamma\Rightarrow \mbox{\boldmath$A$} 
  &  
  \infer{c_{1}: \Gamma,\Delta,\mbox{\boldmath$A$} \Rightarrow C}
   {
   b_{1,1}: \Gamma\Rightarrow B(\bar{n})
   &
   \infer*[D_{r,n}]{c^{\prime}_{1,n}: \Delta, \mbox{\boldmath$A$} , B(\bar{n})\Rightarrow C}{}
  }
 }
 \]
where $b_{1,1}=c_{1}*\langle 2 j_{n}\rangle$ and $c^{\prime}_{1,n}=c_{1}*\langle 2 j_{n}+1\rangle$.

The labeling function is defined by
\begin{equation}\label{eq:labelex5b}
kb_{E}(b_{1,1}*d) = kb_{D_{a}}(b_{1}*d)
\end{equation}
and
\[
kb_{E}(c^{\prime}_{1,n})=kb_{D_{a}}(c_{1,n})
.\]

\item\label{5d}
 The case when $A_{i}\equiv B_{0}\lor B_{1}\in\mbox{\boldmath$A$} $.
 For simplicity suppose $i=1$.
\[
\infer[(L\lor)]
{c_{1}: \Delta,\mbox{\boldmath$A$}  \Rightarrow C}
{
 \infer*[D_{r,0}]{c_{1,0}: \Delta, \mbox{\boldmath$A$} _{1}, B_{0}\Rightarrow C}{}
 &
 \infer*[D_{r,1}]{c_{1,1}: \Delta, \mbox{\boldmath$A$} _{1}, B_{1}\Rightarrow C}{}
 }
 \]
 where $A_{1}\not\in\mbox{\boldmath$A$} _{1}$.
 
 By Inversion Lemma \ref{lem:inversion}.\ref{lem:inversion1} pick an $n=0,1$ such that
$\Gamma\Rightarrow B_{n}$ is provable without changing the naked tree.

Suppose that $n=0$. $E$ is defined as follows.
 \[
 \infer
 {a: \Gamma,\Delta\Rightarrow C}
 {
   \mbox{\boldmath$b$} _{1}: \Gamma\Rightarrow \mbox{\boldmath$A$} _{1}
  &  
  \infer{c_{2}: \Gamma,\Delta,\mbox{\boldmath$A$} _{1}\Rightarrow C}
   {
   b_{1}: \Gamma\Rightarrow B_{0}
   &
   \infer[(Rep)]{c_{1}: \Delta, \mbox{\boldmath$A$} _{1}, B_{0}\Rightarrow C}
   {
   \infer*[D_{r,0}]{c_{1,0}: \Delta, \mbox{\boldmath$A$} _{1}, B_{0}\Rightarrow C}{}
   }
  }
 }
 \]

Note that the new cut inference for $B_{0}$ may be of rank 0.

\item\label{5e}
 The case when $A_{i}\equiv B_{0}\land B_{1} \in\mbox{\boldmath$A$} $ is treated as in the case \ref{5c} for universal quantifier.
 \end{enumerate}
\end{enumerate}

\begin{claim}\label{clm:kb}
The resulting derivation $ce_{1}(D)$ can be labeled enjoying the conditions (\ref{eq:kb1}), (\ref{eq:kb2})
and (\ref{eq:kb3}).
\end{claim}
{\bf Proof}.\hspace{2mm}
Let $D_{a}$ be the trunk ending with the leftmost cut of rank 1 in $D$.
First observe that the labels $\{kb_{E}(b): b\in T(E)\}\subseteq\{kb_{D_{a}}(b):b\in T(D_{a})\}$.
Therefore it suffices to see that $E$, and hence $ce_{2}(\mbox{\boldmath$D$} _{\ell},D_{r})$ enjoys the three conditions
if $D_{a}$ does.

Note that (the naked tree of) $E$ is constructed from $D_{r}$ by appending trees $\mbox{\boldmath$D$} _{\ell}$
only where
 a right cut formula $A_{i}$ has an ancestor which is either a formula of rank 0 or
a principal formula of an initial sequent $\Phi,I(t)\Rightarrow I(t)$.
In the latter case the ancestor has to be the formula $I(t)$ in the antecedent.

$E$ enjoys the first (\ref{eq:kb1}) since $D_{a}$ does the first (\ref{eq:kb1}).
$E$ enjoys the second (\ref{eq:kb2}) since $D_{a}$ does the third (\ref{eq:kb3}).
$E$ enjoys the third (\ref{eq:kb3}) since $D_{a}$ does the third (\ref{eq:kb3}), and
the first (\ref{eq:kb1}).

Let us examine cases.
Consider the case \ref{4} when $\mbox{\boldmath$D$} _{\ell}$ is lifted up.
We have to show
\[
kb_{E}(e)<_{T_{2}}kb_{E}(c^{\prime}_{j,i})
\]
for $e$ such that $b_{j,i}\subseteq e$.
Then $kb_{E}(e)=kb_{E}(a*\langle n_{i}\rangle*\langle r_{k},\ldots, r_{j+1},\ell_{j}\rangle*d)= kb_{D_{a}}(a*\langle r_{k},\ldots, r_{j+1},\ell_{j}\rangle*d)=kb_{D_{a}}(b_{j}*d)$ and $kb_{E}(c^{\prime}_{j,i})=kb_{D_{a}}(c_{j})$.

Since the right cut formula $A_{i}$ has an ancestor in $c_{1,i}:\Delta,\mbox{\boldmath$A$} \Rightarrow C_{i}$,
$kb_{E}(e)<_{T_{2}}kb_{E}(c^{\prime}_{j,i})$ follows from  (\ref{eq:kb3}) for $D_{a}$.

Next consider the case \ref{5b}.
 Although $b_{j,m}, c_{j,m}$ are duplicated for $m=0,1$, and 
 $kb_{E}(c_{j,0}) = kb_{D_{a}}(c_{j})=kb_{E}(c_{j,1})$, 
$kb_{E}(b_{j,1}*d) = kb_{D_{a}}(b_{j}*d)=kb_{E}(b_{j,0}*d)$
by (\ref{eq:labelex}), (\ref{eq:labelex5b})
these are harmless for (\ref{eq:kb3})
since the juncture is a cut of rank 0, $rk(H)=0$.

Finally consider the case \ref{5c}.

By (\ref{eq:labelex5b}) we have
\[
kb_{E}(b_{1,1}*d) = kb_{D_{a}}(b_{1}*d)=kb_{E}(b_{1}*d)
\]
but the right cut formula $A_{i}$ has no ancestor in $b_{1,1}:\Gamma\Rightarrow B(\bar{n})$.
Thus (\ref{eq:kb3}) is enjoyed.

Next for (\ref{eq:kb2}) we have
\[
kb_{D_{a}}(c_{j})<_{T_{2}}kb_{D_{a}}(b_{1}*d)=kb_{E}(b_{1,1}*d) 
\]
by (\ref{eq:kb2}) and (\ref{eq:kb3}) in $D_{a}$.

Finally for (\ref{eq:kb1}) assume $j\neq 1$ and
\[
kb_{D_{a}}(b_{1}*d)=kb_{E}(b_{1,1}*d)\neq kb_{E}(b_{j}*e)= kb_{D_{a}}(b_{j}*e) 
.\]
Then by (\ref{eq:kb1}) in $D_{a}$ we have $b_{j}*e<_{T(D_{a})}b_{1}*d$, and hence
$kb_{D_{a}}(b_{j}*e) <_{T_{2}}kb_{D_{a}}(b_{1}*d)$.
\hspace*{\fill} $\Box$

This ends the construction of the cut-eliminating operator $ce_{1}(D)$.

Finally we show Theorem \ref{lem:quickelim}.
Given a derivation $D_{2}=(T_{2},Seq,Rule, rk,ord,kb)$ of $\Gamma_{0}\Rightarrow C_{0}$ of rank 1, and assume $\Gamma_{0}\subseteq{\cal H}$ and $C_{0}\in {\cal HM}$.
$(T_{2},<_{T_{2}})$ denotes the Kleene-Brouwer ordering on the naked tree $T_{2}$.

Let $KB(D):=KB(J)$ for the leftmost cut inference $J$ of rank 1
if such a $J$ exists.
Otherwise $KB(D)$ denote the largest element in $T_{2}$ with respect to $<_{T_{2}}$,
i.e., the root of $T_{2}$.
Then we see that $KB(D)<_{T_{2}}KB(ce_{1}(D))$ if $D$ contains a cut inference of rank 1.

Suppose as the induction hypothesis that any cut inferences $J$ of rank 1 has been eliminated
for $KB(J)<a$, and let $D$ denote such a derivation.
Also assume that $a$ is a node of a cut inference of rank 1.
Then in $ce_{1}(D)$ the cut inference is eliminated.
This proves the Theorem \ref{lem:quickelim} by induction along the Kleene-Brouwer ordering $<_{T_{2}}$, cf. Lemma \ref{lem:KB}.


\begin{thebibliography}{99}


\bibitem{attic}T. Arai, Some results on cut-elimination, provable well-orderings, induction and reflection, 
Annals of Pure and Applied Logic vol. 95 (1998),  pp. 93-184.

\bibitem{Motohashi} T. Arai, Non-elementary speed-ups in logic calculi,
 Mathematical Logic Quarterly vol. 6(2008),  pp. 629-640.

\bibitem{Buchholz}W. Buchholz, An intuitionistic fixed point theory, Arch. Math. Logic 37(1997), pp. 21-27.

\bibitem{Feferman}S. Feferman, Iterated inductive fixed-point theories:Applications to Hancock's conjecture, in:G. Metakides, ed., Patras Logic Symposion (North-Holland, Amsterdam, 1982), pp. 171-196.

\bibitem{Goodman}N. Goodman, Relativized realizability in intuitionistic arithmetic of all finite types, J. Symb. Logic 43 (1978), pp. 23-44.

\bibitem{Mintsfinite} G.E. Mints, Finite investigations of transfinite derivations, in: Selected Papers in Proof Theory (Bibliopolis, Napoli, 1992), pp. 17-72.

\bibitem{Mintsmono} G. E. Mints, Quick cut-elimination for monotone cuts, in 
Games, logic, and constructive sets(Stanford, CA, 2000), CSLI Lecture Notes, 161, CSLI Publ., Stanford, CA, 2003, pp. 75-83.

\bibitem{Strahm}C. R\"uede and T. Strahm, Intuitionistic fixed point theories for strictly positive operators, 
Math. Log. Quart. 48(2002), pp. 195-202.

\bibitem{Troelstra}A.S. Troelstra, Metamathematical Investigation of Intuitionistic Arithmetic and Analysis. Lecture Notes in Mathematics 344 (Springer, Berlin Heidelberg New York, 1973). 


\end{thebibliography}
\end{document}